\documentstyle[12pt, amsfonts]{amsart}
\begin{document}
\def\R{{\mathbb R}}
\def\Z{{\mathbb Z}}
\def\C{{\mathbb C}}
\newcommand{\trace}{\rm trace}
\newcommand{\Ex}{{\mathbb{E}}}
\newcommand{\Prob}{{\mathbb{P}}}
\newcommand{\E}{{\cal E}}
\newcommand{\F}{{\cal F}}
\newtheorem{df}{Definition}
\newtheorem{theorem}{Theorem}
\newtheorem{lemma}{Lemma}
\newtheorem{pr}{Proposition}
\newtheorem{co}{Corollary}
\def\n{\nu}
\def\sign{\mbox{ sign }}
\def\a{\alpha}
\def\N{{\mathbb N}}
\def\A{{\cal A}}
\def\L{{\cal L}}
\def\X{{\cal X}}
\def\F{{\cal F}}
\def\c{\bar{c}}
\def\v{\nu}
\def\d{\delta}
\def\diam{\mbox{\rm dim}}
\def\vol{\mbox{\rm Vol}}
\def\b{\beta}
\def\t{\theta}
\def\l{\lambda}
\def\e{\varepsilon}
\def\colon{{:}\;}
\def\pf{\noindent {\bf Proof :  \  }}
\def\endpf{ \begin{flushright}
$ \Box $ \\
\end{flushright}}

\title[Estimates for measures of sections]{Estimates for measures of sections of convex bodies}

\author{Alexander Koldobsky}

\address{Department of Mathematics\\ 
University of Missouri\\
Columbia, MO 65211}

\email{koldobskiya@@missouri.edu}

\begin{abstract}  A $\sqrt{n}$ estimate in the hyperplane problem with arbitrary measures 
has recently been proved in \cite{K3}.  In this note we present analogs of this result
for sections of lower dimensions and in the complex case. We deduce these inequalities from 
stability in comparison problems for different generalizations of intersection bodies.

\end{abstract}  
\maketitle

\section{Introduction}
The following inequality has recently been proved in \cite{K3}. Let $K$ be an origin symmetric convex 
body in $\R^n,$ and let $\mu$ be a measure on $K$ with even continuous non-negative density $f$
so that $\mu(B)=\int_B f$ for every Borel subset of $K.$ Then
\begin{equation} \label{sqrtn1}
\mu(K)\ \le\ \sqrt{n} \frac n{n-1} c_n\max_{\xi \in S^{n-1}} 
\mu(K\cap \xi^\bot)\ |K|^{1/n} \ ,
\end{equation}
where  $c_n= \left|B_2^n\right|^{\frac{n-1}n}/ \left|B_2^{n-1}\right| < 1,$
$B_2^n$ is the unit Euclidean ball in $\R^n,$ and $|K|$ stands for volume of proper dimension.
Note that $c_n<1$ for every $n.$

In the case of volume, when $f=1$ everywhere on $K,$ inequality (\ref{sqrtn1}) was proved in
\cite[p. 96]{MP}. Another argument follows from
 \cite[Theorem 8.2.13]{G}; in \cite{G} this argument is attributed to Rolf Schneider.  Also, in the case of volume the constant
$\sqrt{n}$ can be improved to $Cn^{1/4},$ where $C$ is an absolute constant, as shown by Klartag \cite{Kl}
who removed a logarithmic term from an earlier estimate of Bourgain \cite{Bo3}. These results
are much more involved. The question of whether $n^{1/4}$ can also be removed in the case of volume
is the matter of the hyperplane conjecture \cite{Bo1,Bo2, Ba, MP}; see the book \cite{BGVV} for the current state 
of the problem.

In this note we prove analogs of inequality (\ref{sqrtn1}) for sections of lower dimensions and
in the complex case; see Theorems \ref{main1} and \ref{main2},  respectively. As in \cite{K3}, 
the proofs are based on certain stability results for
generalizations of intersection bodies.

\section{Lower dimensional sections}

We need several definitions and facts.
A closed bounded set $K$ in $\R^n$ is called a {\it star body}  if 
every straight line passing through the origin crosses the boundary of $K$ 
at exactly two points different from the origin, the origin is an interior point of $K,$
and the {\it Minkowski functional} 
of $K$ defined by 
$$\|x\|_K = \min\{a\ge 0:\ x\in aK\}$$
is a continuous function on $\R^n.$ 

The {\it radial function} of a star body $K$ is defined by
$$\rho_K(x) = \|x\|_K^{-1}, \qquad x\in \R^n.$$
If $x\in S^{n-1}$ then $\rho_K(x)$ is the radius of $K$ in the
direction of $x.$

If $\mu$ is a measure on $K$ with even continuous density $f$, then 
\begin{equation} \label{polar-measure}
\mu(K) = \int_K f(x)\ dx = \int\limits_{S^{n-1}}\left(\int\limits_0^{\|\theta\|^{-1}_K} r^{n-1} f(r\theta)\ dr\right) d\theta.
\end{equation}
Putting $f=1$, one gets
\begin{equation} \label{polar-volume}
|K|
=\frac{1}{n} \int_{S^{n-1}} \rho_K^n(\theta) d\theta=
\frac{1}{n} \int_{S^{n-1}} \|\theta\|_K^{-n} d\theta.
\end{equation}

For $1\le k \le n-1,$ denote by $Gr_{n-k}$ the Grassmanian of $(n-k)$-dimensional
subspaces of $\R^n.$ The {\it $(n-k)$-dimensional spherical Radon transform} 
$R_{n-k}:C(S^{n-1})\mapsto C(Gr_{n-k})$  
is a linear operator defined by
$$R_{n-k}g (H)=\int_{S^{n-1}\cap H} g(x)\ dx,\quad \forall  H\in Gr_{n-k}$$
for every function $g\in C(S^{n-1}).$

The polar formulas (\ref{polar-measure}) and  (\ref{polar-volume}), applied to sections of $K$, express 
volume in terms of the spherical Radon transform:

$$\mu(K\cap H) = \int_{K\cap H} f =  
\int_{S^{n-1}\cap H} \left(\int_0^{\|\theta\|_K^{-1}} r^{n-k-1}f(r\theta)\ dr \right)d\theta$$
\begin{equation} \label{measure=spherradon}
=R_{n-k}\left(\int_0^{\|\cdot\|_K^{-1}} r^{n-k-1}f(r\ \cdot)\ dr \right)(\xi).
\end{equation}
and
\begin{equation} \label{volume=spherradon}
|K\cap H| = \frac{1}{n-k} \int_{S^{n-1}\cap \xi^\bot} \|\theta\|_K^{-n+k}d\theta =
\frac{1}{n-k} R_{n-k}(\|\cdot\|_K^{-n+k})(\xi).
\end{equation}

The class of intersection bodies was introduced by Lutwak \cite{L} and played a crucial role in the
solution of the Busemann-Petty problem; see \cite{G,K1} for definition and properties. 
A more general class of bodies was introduced by Zhang \cite{Z} 
in connection with the lower dimensional Busemann-Petty problem.
Denote 
$$R_{n-k}\left(C(S^{n-1})\right)=X\subset C(Gr_{n-k}).$$
Let $M^+(X)$ be the space of linear positive continuous functionals on $X$, i.e. for every
$\nu\in M^+(X)$ and non-negative function $f\in X$, we have $\nu(f)\geq0$.

An origin-symmetric star body $K$ in $\R^n$ is called a {\it generalized $k$-intersection
body} if there exists a functional $\nu\in M^+(X)$ so that for every $g\in C(S^{n-1})$,
\begin{equation} \label{defintbody}
\int_{S^{n-1}} \|x\|_K^{-k} g(x)\ dx=\nu(R_{n-k}g).
\end{equation}

When $k=1$ we get the class of intersection bodies.
It was proved by Grinberg and Zhang \cite[Lemma 6.1]{GZ} that every intersection body in $\R^n$
is a generalized $k$-intersection
body for every $k<n.$ More generally, as proved later by Milman \cite{Mi}, if  $m$ divides $k$, then every
generalized $m$-intersection body is a generalized $k$-intersection body.  

We need the following stability result for generalized $k$-intersection bodies.  
\begin{theorem}\label{stab1}
Suppose that $1\le k \le n-1,$ $K$ is a generalized $k$-intersection body in $\R^n,$  $f$
is an even continuous function on $K,$ $f\ge 1$ everywhere on $K,$ and $\e>0.$ If
\begin{equation}\label{comp1}
\int_{K\cap H} f \ \le\ |K\cap H| +\e,\qquad \forall H\in Gr_{n-k},
\end{equation}
then
\begin{equation}\label{comp2}
\int_K f\ \le\ |K| + \frac {n}{n-k}\ c_{n,k}\ |K|^{k/n}\e,
\end{equation}
where $c_{n,k}= |B_2^n|^{\frac {n-k}n}/|B_2^{n-k}| < 1.$
\end{theorem}

\pf  Use polar formulas (\ref{measure=spherradon}) and (\ref{volume=spherradon}) to write
the condition (\ref{comp1}) in terms of the $(n-k)$-dimensional spherical Radon transform: for all $H\in Gr_{n-k}$
$$R_{n-k}\left(\int_0^{\|\cdot\|_K^{-1}} r^{n-k-1}f(r\ \cdot)\ dr \right)(H) \le \frac{1}{n-k} R_{n-k} \left(\|\cdot\|_K^{-n+k} \right)(H) + \e.$$
Let $\nu$ be the functional corresponding to $K$ by (\ref{defintbody}), apply $\nu$ to both sides of the latter inequality
(the direction of the inequality is preserved because $\nu$ is a positive functional) and use (\ref{defintbody}).
We get
$$\int_{S^{n-1}} \|\theta\|_K^{-k} \left(\int_0^{\|\theta\|_K^{-1}} r^{n-k-1}f(r\theta)\ dr \right)d\theta $$
\begin{equation}\label{eq00}
\le  \frac{1}{n-k} \int_{S^{n-1}} \|\theta\|_K^{-n}\ d\theta + \e \nu(1).
\end{equation}
Split the integral in the left-hand side into two integrals and then use $f\ge 1$ as follows: 
$$\int_{S^{n-1}}  \left(\int_0^{\|\theta\|_K^{-1}} r^{n-1}f(r\theta)\ dr \right)d\theta$$ 
$$+ \int_{S^{n-1}} \left(\int_0^{\|\theta\|_K^{-1}} (\|\theta\|_K^{-k} - r^k)  r^{n-k-1}f(r\theta)\ dr \right)d\theta$$
$$\ge \int_K f + \int_{S^{n-1}} \left(\int_0^{\|\theta\|_K^{-1}} (\|\theta\|_K^{-k} - r^k)  r^{n-k_1}\ dr \right)d\theta$$
\begin{equation}\label{eq11}
= \int_K f + \frac1{n-k} |K|.
\end{equation}

Now estimate $\nu(1)$ by first writing $1= R_{n-k}1/|S^{n-k-1}|$ and then using definition (\ref{defintbody}), H\"older's  
inequality and $|S^{n-1}|=n|B_2^n|$:
$$ \nu(1)= \frac 1{\left|S^{n-k-1}\right|} \nu(R_{n-k}1)
=\frac 1{\left| S^{n-k-1} \right| } \int_{S^{n-1}} \|\theta\|_K^{-k}\ d\theta $$
$$ \le  \frac 1{\left|S^{n-k-1}\right|} \left|S^{n-1}\right|^{\frac{n-k}n} \left(\int_{S^{n-1}} \|\theta\|_K^{-n}\ d\theta\right)^{\frac kn}$$
\begin{equation}\label{eq22}
=  \frac{\e}{\left|S^{n-k-1}\right|} \left|S^{n-1}\right|^{\frac{n-k}n} n^{k/n}|K|^{k/n}= \frac n{n-k} c_{n,k} |K|^{k/n}.
\end{equation}
Combining (\ref{eq00}), (\ref{eq11}) and (\ref{eq22})  we get
$$\int_K f + \frac 1{n-k} |K| \le \frac n{n-k} |K| + \frac n{n-k} c_{n,k} |K|^{k/n} \e. \qed$$
\bigbreak
It was proved in \cite{KM} (generalizing the result for $k=1$ from \cite{K2})
that if $L$ is a generalized $k$-intersection body
and $\mu$ is a measure with even continuous density, then
$$\mu(L)\,\leq\,\frac{n}{n-k}c_{n,k}\max_{H\in Gr_{n-k}} \mu(L\cap H) \ |L|^{k/n}.$$
We show now that it is possible to extend this inequality to arbitrary origin-symmetric convex bodies
in $\R^n$ at the expense of an extra constant $n^{k/2}.$
\begin{theorem}\label{main1} Suppose that $L$ is an origin-symmetric convex body in $\R^n,$ and 
$\mu$ is a measure with even continuous non-negative density $g$ on $L.$ Then
\begin{equation} \label{sqrtn}
\mu(L)\ \le\  n^{k/2} \frac n{n-k} c_{n,k} \max_{H\in Gr_{n-k}}  \mu(L\cap H)\ |L|^{k/n}.
\end{equation}
\end{theorem}

\pf By John's theorem \cite{J}, there exists an origin-symmetric ellipsoid $K$ such that
$$\frac 1{\sqrt{n}} K \subset L \subset K.$$ 
The ellipsoid $K$ is an intersection body (\cite[Corollary 8.1.7]{G}),
and every intersection body is a generalized $k$-intersection body for every $k$ (\cite[Lemma 6.1]{GZ}).
Let $f= \chi_K + g \chi_L,$ where $\chi_K,\ \chi_L$ are the indicator functions of $K$ and $L,$ 
then $f\ge 1$ everywhere on $K.$ Put 
$$\e=\max_{H\in Gr_{n-k}} \left(\int_{K\cap H} f - |K\cap H| \right)= \max_{H\in Gr_{n-k}} \int_{L\cap H} g.$$
Now we can apply Theorem \ref{stab1} to $f,K,\e$ (the function $f$ is not necessarily continuous on $K,$ 
but the result holds by a simple approximation argument). We get
$$\mu(L)= \int_L g = \int_K f  -\ |K|$$
$$ \le \frac n{n-k} c_{n,k}  |K|^{k/n}\max_{H\in Gr_{n-k}} \int_{L\cap H} g$$
$$ \le n^{k/2} \ \frac n{n-k} c_{n,k} |L|^{k/n}\max_{H\in Gr_{n-k}} \mu(L\cap H),$$
because $K\subset \sqrt{n}L,$ so $|K|\le n^{n/2}  |L|.$ \qed 

\section{The complex case}
Origin symmetric convex bodies in $\C^n$ are the unit balls of norms on $\C^n.$
We denote by $\|\cdot\|_K$
the norm corresponding to the body $K:$
$$K=\{z\in \C^n:\ \|z\|_K\le 1\}.$$
In order to define volume, we identify $\C^n$ with $\R^{2n}$ using the standard mapping
$$\xi = (\xi_1,...,\xi_n)=(\xi_{11}+i\xi_{12},...,\xi_{n1}+i\xi_{n2})
 \mapsto  (\xi_{11},\xi_{12},...,\xi_{n1},\xi_{n2}).$$
 Since norms on $\C^n$ satisfy the equality
$$\|\lambda z\| = |\lambda|\|z\|,\quad \forall z\in \C^n,\  \forall\lambda \in \C,$$
origin symmetric complex convex bodies correspond to those origin symmetric convex bodies
$K$  in $\R^{2n}$ that are invariant
 with respect to any coordinate-wise two-dimensional rotation, namely for each $\theta\in [0,2\pi]$
 and each $\xi= (\xi_{11},\xi_{12},...,\xi_{n1},\xi_{n2})\in \R^{2n}$
  \begin{equation} \label{rotation}
  \|\xi\|_K =
 \|R_\theta(\xi_{11},\xi_{12}),...,R_\theta(\xi_{n1},\xi_{n2})\|_K,
 \end{equation}
 where $R_\theta$ stands for  the counterclockwise rotation of $\R^2$ by the angle
 $\theta$ with respect to the origin. We shall say that $K$ is a {\it complex convex body
 in $\R^{2n}$} if $K$ is a convex body and satisfies equations (\ref{rotation}). Similarly,
 complex star bodies are $R_\theta$-invariant star bodies in $\R^{2n}.$
 \medbreak
 For $\xi\in \C^n,
|\xi|=1,$ denote by
$$H_\xi = \{ z\in \C^n:\ (z,\xi)=\sum_{k=1}^n z_k\overline{\xi_k} =0\}$$
the complex hyperplane through the origin, perpendicular to $\xi.$
 Under the standard mapping from $\C^n$ to $\R^{2n}$ the hyperplane $H_\xi$ 
 turns into a $(2n-2)$-dimensional subspace of $\R^{2n}.$ 
   \medbreak
Denote by $C_c(S^{2n-1})$ the space of $R_\theta$-invariant continuous functions, i.e.
  continuous real-valued functions $f$ on the unit sphere $S^{2n-1}$ in $\R^{2n}$ satisfying 
  $f(\xi)=f(R_\theta(\xi))$ for all $\xi\in S^{2n-1}$ and all $\theta\in [0,2\pi].$ The {\it complex spherical
  Radon transform} is an operator ${\cal{R}}_c: C_c(S^{2n-1})\to C_c(S^{2n-1})$ defined by
  $${\cal{R}}_cf(\xi) = \int_{S^{2n-1}\cap H_\xi} f(x) dx.$$
 
  We say that a finite Borel measure $\mu$ on
$S^{2n-1}$ is $R_\theta$-invariant if for any continuous function $f$ on $S^{2n-1}$ and any $\theta\in [0,2\pi]$,
$$\int_{S^{2n-1}} f(x) d\mu(x) = \int_{S^{2n-1}} f(R_\theta x) d\mu(x).$$
The complex spherical Radon transform of an $R_\theta$-invariant measure $\mu$ is defined
as a functional ${\cal{R}}_c\mu$ on the space $C_c(S^{2n-1})$ acting by 
$$ \left({\cal{R}}_c\mu, f \right) = \int_{S^{2n-1}} {\cal{R}}_cf(x) d\mu(x).$$

Complex intersection bodies were introduced and studied in \cite{KPZ}. 
An origin symmetric complex star body $K$ in $\R^{2n}$ is called a {\it complex intersection body} if there
exists a finite Borel $R_\theta$-invariant measure $\mu$ on $S^{2n-1}$ so that
$\|\cdot\|_K^{-2}$ and ${\cal{R}}_c\mu$ are equal as functionals on $C_c(S^{2n-1}),$ i.e.
for any $f\in C_c(S^{2n-1})$
\begin{equation}\label{defcompint}
\int_{S^{2n-1}} \|x\|_K^{-2} f (x)\ dx = \int_{S^{2n-1}} {\cal{R}}_c f(\theta) d\mu(\theta).
\end{equation}
\bigbreak
\begin{theorem} \label{stab3} Suppose that $K$ is a complex intersection body in $\R^{2n},$
$f$ is an even continuous $R_\theta$-invariant function on $K,$ $f\ge 1$ everywhere on $K,$
and $\e>0.$ If
\begin{equation}\label{comp3}
\int_{K\cap H_\xi} f \ \le\ |K\cap H_\xi| +\e,\qquad \forall \xi\in S^{2n-1},
\end{equation}
then
\begin{equation}\label{comp4}
\int_K f\ \le\ |K| + \frac {n}{n-1}\ d_n\ |K|^{1/n}\e,
\end{equation}
where $d_n= |B_2^{2n}|^{\frac {n-1}n}/|B_2^{2n-2}| < 1.$
\end{theorem}

\pf  Use the polar formulas (\ref{measure=spherradon}) and (\ref{volume=spherradon}) to write
the condition (\ref{comp3}) in terms of the complex spherical Radon transform: for all $\xi\in S^{2n-1}$
$$ {\cal{R}}_c \left(\int_0^{\|\cdot\|_K^{-1}} r^{2n-3}f(r\ \cdot)\ dr \right)(\xi) \le 
\frac{1}{2n-2}  {\cal{R}}_c \left(\|\cdot\|_K^{-2n+2}\right)(\xi) + \e.$$
Let $\mu$ be the measure on $S^{2n-1}$ corresponding to $K$ by (\ref{defcompint}). Integrate
the latter inequality over $S^{2n-1}$ with the measure $\mu$ and use (\ref{defcompint}):
$$\int_{S^{2n-1}} \|\theta\|_K^{-2} \left(\int_0^{\|\theta\|_K^{-1}} r^{2n-3}f(r\theta)\ dr \right)d\theta $$
$$\le  \frac{1}{2n-2} \int_{S^{2n-1}} \|\theta\|_K^{-2n}\ d\theta + \e \int_{S^{2n-1}} d\mu(\xi)$$
\begin{equation} \label{eq44}
= \frac n{n-1} |K| +  \e \int_{S^{2n-1}} d\mu(\xi) .
\end{equation}
Recall (\ref{polar-measure}), (\ref{polar-volume}) and the assumption that $f\ge 1.$ We estimate the integral in the left-hand side 
of (\ref{eq44}) as follows:
$$\int_{S^{2n-1}} \|\theta\|_K^{-2} \left(\int_0^{\|\theta\|_K^{-1}} r^{2n-3}f(r\theta)\ dr \right)d\theta $$
$$= \int_{S^{2n-1}}  \left(\int_0^{\|\theta\|_K^{-1}} r^{2n-1}f(r\theta)\ dr \right)d\theta$$ 
$$+ \int_{S^{2n-1}} \left(\int_0^{\|\theta\|_K^{-1}} (\|\theta\|_K^{-2} - r^2)  r^{2n-3}f(r\theta)\ dr \right)d\theta$$
$$\ge \int_K f + \int_{S^{2n-1}} \left(\int_0^{\|\theta\|_K^{-1}} (\|\theta\|_K^{-2} - r^2)  r^{2n-3}\ dr \right)d\theta$$
\begin{equation} \label{eq55}
=\int_K f + \frac 1{2(n-1)n} \int_{S^{2n-1}} \|\theta\|_K^{-2n}\ d\theta = \int_K f + \frac1{n-1} |K|.
\end{equation}

Let us estimate the second term in the right-hand side of (\ref{eq44}) by adding the complex spherical 
Radon transform of  the unit constant
function under the integral (${\cal{R}}_c1(\xi)=\left|S^{2n-3}\right|$ for every $\xi \in S^{2n-1}$), 
using again (\ref{defcompint}) and then applying H\"older's  inequality:
$$\e \int_{S^{2n-1}} d\mu(\xi) = \frac{\e}{\left|S^{2n-3}\right|} \int_{S^{2n-1}} {\cal{R}}_c1(\xi)\ d\mu(\xi)$$
$$=\frac{\e}{\left| S^{2n-3} \right| } \int_{S^{2n-1}} \|\theta\|_K^{-2}\ d\theta $$
$$ \le  \frac{\e}{\left|S^{2n-3}\right|} \left|S^{2n-1}\right|^{\frac{n-1}n} \left(\int_{S^{2n-1}} \|\theta\|_K^{-2n}\ d\theta\right)^{\frac1n}$$
\begin{equation}\label{eq66}
=  \frac{\e}{\left|S^{2n-3}\right|} \left|S^{2n-1}\right|^{\frac{n-1}n} (2n)^{1/n}|K|^{1/n}= \frac n{n-1} d_n |K|^{1/n} \e.
\end{equation}
In the last step we used $|S^{2n-1}|=2n|B_2^{2n}|.$ Combining (\ref{eq44}),(\ref{eq55}),(\ref{eq66}) we get
$$\int_K f + \frac 1{n-1} |K| \le \frac n{n-1} |K| + \frac n{n-1} d_n |K|^{1/n} \e. \qed$$

It was proved in \cite{KPZ} that if $K$ is a complex intersection body in $\R^{2n}$ and $\gamma$ 
is an arbitrary measure on $\R^{2n}$ with even continuous density, then
$$\gamma(K) \le \frac n{n-1} d_n \max_{\xi\in S^{2n-1}} \gamma(K\cap H_\xi)\ |K|^{\frac1n}.$$
We remove the condition that $K$ is a complex intersection body at the expense of an extra
constant.

We use a result from \cite[Theorem 2]{KPZ} that a complex star body is a complex intersection body
if and only if $\|\cdot\|_K^{-2}$ is a positive definite distribution, i.e. its Fourier transform
in the sense of distributions assumes non-negative values on non-negative test functions.
We refer the reader to \cite{K1,KPZ} for details. 

\begin{theorem}  \label{main2} Suppose that $L$ is an origin-symmetric complex convex body in $\R^{2n}$ and $\gamma$ 
is an arbitrary measure on $\R^{2n}$ with even continuous density $g$, then
$$\gamma(L) \le 2n \frac {n}{n-1} d_n \max_{\xi\in S^{2n-1}} \gamma(L\cap H_\xi)\ |L|^{\frac1n}.$$
\end{theorem}

\pf By John's theorem \cite{J}, there exists an origin-symmetric ellipsoid $K$ such that
$$\frac 1{\sqrt{2n}} K \subset L \subset K.$$ 
Construct a new body $K_c$ by
$$\|x\|_{K_c}^{-2}= \frac 1{2\pi} \int_0^{2\pi} \|R_\theta x\|_K^{-2} d\theta.$$
Clearly, $K_c$ is $R_\theta$-invariant, so it is a complex star body. For every $\theta\in [0,2\pi]$ the 
distribution $\|R_\theta x\|_K^{-2}$ is positive definite, because this is 
a linear transformation of the Euclidean norm. So $\|x\|_{K_c}^{-2}$ is
also a positive definite distribution, and, by \cite[Theorem 2]{KPZ}, 
$K_c$ is a complex intersection body. Since $\frac 1{\sqrt{2n}} K\subset L \subset K$
and $L$ is $R_\theta$-invariant as a complex convex body, we have 
$$\frac 1{\sqrt{2n}} R_\theta K \subset L \subset R_\theta K, \quad \forall \theta\in [0,2\pi],$$
so
$$\frac 1{\sqrt{2n}} K_c\subset L \subset K_c.$$
Let $f= \chi_{K_c} + g \chi_L,$ where $\chi_{K_c},\ \chi_L$ are the indicator functions of $K_c$ and $L.$ 
Clearly, $f$ is $R_\theta$-invariant and  $f\ge 1$ everywhere on $K.$ Put 
$$\e=\max_{\xi\in S^{2n-1}} \left(\int_{K_c\cap H_\xi} f - |K_c\cap H_\xi| \right)= \max_{\xi\in S^{2n-1}} \int_{L\cap H_\xi} g$$
and apply Theorem \ref{stab3} to $f,K_c,\e$ (the function $f$ is not necessarily continuous on $K_c,$ 
but the result holds by a simple approximation argument). We get
$$\mu(L)= \int_L g = \int_{K_c} f  -\ |K_c|$$
$$ \le \frac n{n-1} d_n |K_c|^{1/n}\max_{\xi\in S^{2n-1}} \int_{L\cap H_\xi} g$$
$$ \le 2n \ \frac n{n-1} d_n |L|^{1/n}\max_{\xi\in S^{2n-1}} \mu(L\cap H_\xi),$$
because $|K_c|^{1/n}\le 2n\ |L|^{1/n}.$ \qed 
\bigbreak
Theorem \ref{main2} shows that if bodies have additional symmetries then maximum
in the slicing inequality can be taken over a rather small set of subspaces. 

\bigbreak
{\bf Acknowledgement.} I wish to thank the US National Science Foundation for support through 
grant DMS-1265155.

\end{document}